\documentclass{amsart}

\usepackage[english]{babel}

\usepackage[normalem]{ulem} 
\usepackage[super]{nth} 
\usepackage{setspace}
\usepackage{amsmath, amsthm,amssymb}
\usepackage{float}
\usepackage{bm}
\usepackage{graphicx}
\usepackage{pgfplots}
\usepackage{pgfplotstable}
\pgfplotsset{compat = newest}

\begin{document}

\title{On the Leibnitz Rule for Differentiating Under the Integral Sign}

\author{Jean-Luc Boulnois}

\address{Babson College, Babson Park, Wellesley, Massachusetts 02457, jlboulnois@msn.com}

\begin{abstract}
This  Note revisits the Leibnitz integral calculus method based on differentiation under the integral sign with respect to a  parameter either already existing or introduced \textit{ad hoc}. Through several cases exemplifying the method, it is shown that this approach, applicable to regular and, under certain conditions, to improper integrals as well, results in a first order differential equation whose solution is usually straightforward. 
\end{abstract}

\keywords{Differentiation \and Integrals}

\subjclass[2020]{33B99}
\maketitle


\section{Introduction}

A powerful integral calculus method which takes advantage of differentiation under the integral sign is based on  Leibnitz integral rule which requires that a function $f(x, \alpha)$ of the dependent variable $x$ to be integrated over some interval together with  its derivative be continuous over the integration interval to allow for the interchange of the integration and differentiation operators \cite{hildebrand}. 
\newline

This rule, applicable to definite or indefinite integrals,  also applies to improper integrals where one or both limits are $\infty$. In this particular case, a condition for the differentiation under the integral sign with respect to $\alpha$ to be valid, is that first there exists a function $M(\alpha)$ independent of $x$ such that $|\partial{f} / \partial{ \alpha}| \leqslant M(\alpha) $, and second that the integral over $\alpha$ converges over the integration interval  \cite{hildebrand}. 
\newline

The Leibnitz integral rule can be generalized to multidimensional integrals \cite{flanders}.
\newline

At some point, mathematicians realized that once differentiation had been performed, the Leibnitz rule could be extended to calculate the original integral. The method was popularized by R. P. Feynman as "the Feynman trick" \cite{feynman}, yet it has long been taught in most Advanced Calculus courses in France \cite{pharabod}. 
\newline

The integration method proceeds in several consecutive steps:

\begin{enumerate}
\item First consider the original integral of $f(x, \alpha)$ over some integration interval as a function of $\alpha$; the parameter $\alpha$ is sometimes already present in the integrand, whereas at other times it is judiciously inserted \textit{"ad hoc"} in such a way that derivation with respect to  $\alpha$  yields a second integral which is simpler to evaluate as a function of $x$. 
\item Evaluate the original $x$-integral for a particular convenient value of $\alpha$ (e.g. $\alpha = 0$, $\alpha = 1$, $\alpha = -1$, etc.).
\item Establish the validity of the the method with regard to the Leibnitz  rule for integrating the function $f(x, \alpha)$ as a function of $x$ over the interval.
\item As a function of $x$ evaluate the integral resulting from the differentiation with respect to $\alpha$: this results in a first-order ordinary differential equation in $\alpha$.
\item Subsequently solve  the differential equation yielding a solution to the original integral valid for any value of $\alpha$ up to a constant of integration.
\item Lastly calculate the integration constant  by setting  $\alpha$  to the particular numerical value chosen above: this provides the final solution to the original integral.
\end{enumerate}

In the following, four examples with increasing difficulty are presented to demonstrate the method's versatility while carefully implementing the Leibnitz rule. 
\newline

Also, repeated use will be made of the well-known improper Gauss  integral with parameter $\alpha \in \mathbb{C}$: 

\begin{equation} \label{EQ_1_} 
       \int _{0 }^{\infty } e^{-\alpha x^{2}} dx   = \frac{1}{2} \sqrt{\frac{\pi}{\alpha}}
\end{equation} 
\newline

\section{Example 1}

In order to understand the principle of the method described above, let us first consider the straightforward case of a simple improper integral which  obviously can also be solved by standard integration methods

\begin{equation} \label{EQ_2_}
	I = \int_{0}^{\infty} \ln(1 + x^{2}) \frac{dx}{x^{2}}
\end{equation}

It is noted that the integrand is continuous and bounded since $ \frac {\ln(1+x^{2})}{x^{2}} \leqslant 1$, so the Leibnitz method is applicable. We now freely insert a real positive \textit{ad hoc} parameter $\alpha$ and consider the $\alpha$-dependent auxiliary integral

\begin{equation} \label{EQ_3_}
	J(\alpha) = \int_{0}^{\infty} \ln(1 + \alpha x^{2}) \frac{dx}{x^{2}}
\end{equation}

Obviously $J(\alpha)$ satisfies $ J(1) = I $. Also notice that when $\alpha = 0$ then $J(\alpha) = 0$. Upon differentiating $J(\alpha)$ with respect to $\alpha$, the factor $x^{2}$ in the denominator is eliminated yielding

\begin{equation} \label{EQ_4_}
	\frac {dJ}{d\alpha} = \int_{0}^{\infty} \frac{dx}{1 + \alpha x^{2}} = \frac{1}{\sqrt{\alpha}} \int_{0}^{\infty} \frac{dy}{1 + y^{2}}
\end{equation}

The second integral is obtained by changing variables and letting $y = \sqrt{\alpha} x$. This integral is trivial: its solution is equal to $\arctan(y)$, which for $y = \infty$ is $\pi/2$. Consequently \eqref{EQ_4_} reduces to a simple first-order ordinary differential equation in $\alpha$
	
\begin{equation} \label{EQ_5_}
	\frac {dJ}{d\alpha} = \frac{\pi}{2\sqrt{\alpha}} 
\end{equation}

whose solution is

\begin{equation} \label{EQ_6_}
	J(\alpha) = \pi \sqrt{\alpha} + C
\end{equation}

Since $J(0) = 0$, the integration constant $C$ vanishes and the final solution is obtained by setting $\alpha = 1$ in \eqref{EQ_6_}, yielding

\begin{equation} \label{EQ_7_}
	I = \int_{0}^{\infty} \ln(1 + x^{2}) \frac{dx}{x^{2}} = \pi
\end{equation}

Although this case can evidently be solved by other standard methods, the simple steps of the method described in the introduction and applied here clearly illustrate the value of this approach.

\section{Example 2}

Consider the following often cited example taken from F. Woods whose book \cite{woods}  is mentioned by Feynman \cite{cantor}. This very challenging definite integral depends on an arbitrary parameter $\alpha$ already present in the integrand

\begin{equation} \label{EQ_8_} 
        I(\alpha) =\int _{0}^{\pi } \ln \bigl( 1 - 2 \alpha \cos x + \alpha^{2}\bigr) dx   \quad \text{ for }  |\alpha| \geqslant 1
\end{equation}

It is observed that when $\alpha \gg 1$ the integrand is approximately $\ln(\alpha^{2})$ which provides an upper bound to the integral as $I(\alpha) \leqslant 2 \pi \ln|\alpha|$. Since the integral is bounded, and since the integrand and its derivative are  continuous functions of $x$ over the interval, the Leibnitz rule can be applied. Also, when $\alpha = 1$ the integral vanishes entirely and $I(1) = 0$.
\newline

One approach among others is to introduce the representation of $\cos x$ in terms of complex conjugate exponential functions

\begin{equation} \label{EQ_9_} 
        1 - 2 \alpha \cos x + \alpha^{2} = (\alpha - e^{ix})(\alpha - e^{-ix})
\end{equation}

With $c.c$ as the complex conjugate, the original integral \eqref{EQ_8_} can now be written

\begin{equation} \label{EQ_10_} 
        I(\alpha) =\int _{0}^{\pi} \ln (\alpha - e^{-ix}) dx + c.c.
\end{equation}

Upon performing the differentiation of $I(\alpha)$ with respect to the parameter $\alpha$  one derives

\begin{equation} \label{EQ_11_} 
       \frac{dI}{d\alpha} =  \int _{0}^{\pi} \frac{dx}{\alpha - e^{-ix}}  + c.c. 
\end{equation} 

This can be rearranged as

\begin{equation} \label{EQ_12_} 
       \frac{dI}{d\alpha} =   \frac{1}{\alpha} \int _{0}^{\pi } \Bigl( 1+ \frac{1}{\alpha e^{ix} -1} \Bigr) dx + c.c
\end{equation} 

Therefore \eqref{EQ_12_} reduces to  

\begin{equation} \label{EQ_13_} 
       \frac{dI}{d\alpha} = \frac{2 \pi}{\alpha} + \int _{0}^{\pi }  \frac{e^{-ix}}{\alpha - e^{-ix}} dx + c.c.
\end{equation} 

Calculating  the integral is straightforward, yielding a pure imaginary function of $\alpha$

\begin{equation} \label{EQ_14_} 
	\int _{0}^{\pi }  \frac{e^{-ix}}{\alpha - e^{-ix}} dx = \frac{1}{i} ln \Bigl( \frac{\alpha + 1}{\alpha -1} \Bigr)
\end{equation}

Being a pure imaginary number this integral cancels out with its complex conjugate in \eqref{EQ_13_}, yielding a simple differential equation in $\alpha$

\begin{equation} \label{EQ_15_} 
       \frac{dI}{d\alpha} =  \frac{2\pi}{\alpha} 
\end{equation} 

This is recognized to be the derivative of the solution $\ln(\alpha)$: it includes a constant of integration which vanishes upon setting $\alpha = 1$ as observed above. The final result, identical to the upper limit for $I(\alpha)$, is thus

\begin{equation} \label{EQ_16_} 
	I(\alpha) = 2 \pi \ln|\alpha| \quad \text{ where }  |\alpha| \geqslant 1
\end{equation}

\section{Example 3}

In order to further appreciate the versatility of the method, consider the following non-trivial improper integral

\begin{equation} \label{EQ_17_} 
        I =\int _{0}^{\infty }\frac{e^{-x^{2}} \sin(x^{2}) }{x^{2}} dx  
\end{equation} 

Clearly the integrand in \eqref{EQ_17_} is continuous and, upon setting the exponential to 1, it is bounded with the integral converging absolutely. Since this integral does not contain a "built-in" parameter, the process consists in inserting a real positive \textit{ad hoc} parameter $"\beta"$ either in the exponential or the sine function such that a derivative with respect to $\beta$ eliminates the factor $x^2$ in the denominator. 
\newline

Upon choosing the sine function, one defines the auxiliary integral  $I(\beta)$, function of the parameter $\beta$, for which we seek the solution when $\beta = 1$

\begin{equation} \label{EQ_18_} 
        I(\beta) =\int _{0 }^{\infty }\frac{e^{-x^{2}} \sin(\beta x^{2}) }{x^{2}} dx  
\end{equation} 

The validity of the Leibnitz rule applied to this improper integral is established by observing that, per Eq. \eqref{EQ_1_},  since $|\cos( \beta x^{2})| \leqslant 1 $,  the derivative $dI/d\beta$ below  is bounded by $\sqrt(\pi) /2 $,  independent of $x$ 

\begin{equation} \label{EQ_19_} 
       \frac{dI}{d\beta}=\int _{0 }^{\infty } e^{-x^{2}} \cos(\beta x^{2}) dx   
\end{equation} 

Upon further expressing the cosine function in terms of standard complex exponential functions, with $c. c.$ as the complex conjugate one derives

\begin{equation} \label{EQ_20_} 
       \frac{dI}{d\beta}=\frac{1}{2} \int _{0 }^{\infty } e^{- (1+i \beta) x^{2}} dx   + c. c.
\end{equation} 

Upon integrating over $x$, equation \eqref{EQ_20_} together with the Gaussian integral \eqref{EQ_1_} yields a differential equation in $\beta$ (similar in form to \eqref{EQ_5_})

\begin{equation} \label{EQ_21_} 
       \frac{dI}{d\beta}=\frac{\sqrt{\pi}}{4} \frac{1}{\sqrt{1+i \beta}}   + c. c.
\end{equation} 

Together with its complex conjugate this equation is readily integrated, yielding
a real function for $I(\beta)$ 

\begin{equation} \label{EQ_22_} 
       I(\beta) = i \frac{\sqrt{\pi}}{2} (\sqrt{1-i \beta} - \sqrt{1+i \beta}) + C
\end{equation} 

It is obvious that  setting $\beta = 0$ in \eqref{EQ_18_} implies $I(\beta) = 0$, therefore the above integration constant $C$ must vanish. Upon further expressing the above square roots of conjugate complex numbers in terms of complex exponential functions of a half-angle $\frac{\theta}{2}$ (with $k = 0, 1$), one derives

\begin{equation} \label{EQ_23_} 
       I(\beta) = i \frac{\sqrt{\pi}}{2} (1+\beta^{2})^{\frac{1}{4}}(e^{-i(\frac{\theta}{2}+k\pi)} - e^{i(\frac{\theta}{2}+k\pi)})   \quad \text {with}  \tan\theta = \beta
\end{equation}

This solution is directly expressed in terms of the sine function as

\begin{equation} \label{EQ_24_} 
       I(\beta) = \sqrt{\pi} (1+\beta^{2})^{\frac{1}{4}}  \sin {(\frac {\theta}{2}+k\pi)}  \quad \text {with  k = 0, 1}
\end{equation}

Since \eqref{EQ_23_} is clearly a positive function, only the solution $k = 0$ is retained; expressing $sin(\theta/2)$ in terms of $\beta$ yields the solution

\begin{equation} \label{EQ_25_} 
       I(\beta) = \sqrt{\frac{\pi}{2}}\sqrt{\sqrt{ 1 + \beta^{2}} -1}
\end{equation}

The solution of the  original integral \eqref{EQ_17_} is thus recovered by setting $ \beta = 1$ in \eqref{EQ_25_}. Alternatively,  in equation \eqref{EQ_23_} it is observed that setting $\beta = 1$ yields $\theta = \pi/4$ and thus $\theta /2 = \pi/8$ in \eqref{EQ_24_}, resulting  in the final solution
 
\begin{equation} \label{EQ_26_} 
       I = \sqrt{\frac{\pi}{2}}\sqrt{\sqrt{2} -1}
\end{equation}

Note that if, instead of inserting the parameter $\beta$ in \eqref{EQ_18_}, we had inserted the \textit{ad hoc} parameter $\alpha \geqslant 0$ in the exponential factor, and upon following the same reasoning as above, we would have obtained the solution

\begin{equation} \label{EQ_27_} 
        I(\alpha) =\int _{0}^{\infty }\frac{e^{- \alpha x^{2}} \sin(x^{2}) }{x^{2}} dx  = \sqrt{\frac{\pi}{2}}\sqrt{\sqrt{\alpha^{2}+1} -\alpha}
\end{equation} 

Setting $\alpha = 1$ obviously  yields  solution \eqref{EQ_26_} obtained above. But interestingly, setting $\alpha = 0$ while extending the limits of the integral to $(-\infty; +\infty)$ yields the well known result

\begin{equation} \label{EQ_28_} 
        \int _{-\infty }^{\infty }\frac{\sin( x^{2}) }{x^{2}}  dx  = \sqrt{2 \pi}
\end{equation} 
\newline

\section{Example 4}

The power of the method is best exemplified by considering the following rather difficult definite integral

\begin{equation} \label{EQ_29_} 
        I(a) =\int _{-a}^{+a } \frac{\ln ( 1 + x)}{\sqrt (a^2 - x^2)} dx   \quad \text{ with }  a \in \mathbb{R}  \quad \text{and }  a \geqslant  0  
\end{equation}

For $|x| < a$ this integral is  bounded and continuous. Upon freely inserting an \textit{ad hoc} real positive parameter $"b"$  in the logarithm factor, \eqref{EQ_29_} becomes

\begin{equation} \label{EQ_30_} 
        I(a, b) =\int _{-a}^{+a } \frac{\ln ( 1 + b x)}{\sqrt (a^2 - x^2)} dx   
\end{equation}

We then perform a straightforward elimination of the square root factor in the denominator by changing variables according to

\begin{equation} \label{EQ_31_} 
        x = a sin{\phi}   \quad \text{ which implies  }  \phi = \pm \frac{\pi}{2}
\end{equation}

Upon introducing the parameter $\alpha = ab$ and renaming $I(a, b)$ as $J(\alpha)$, we are then seeking the solution of \eqref{EQ_30_} when letting $b = 1$.  Then $J(\alpha)$  becomes

\begin{equation} \label{EQ_32_} 
        J(\alpha) =\int _{-\frac{\pi}{2}}^{\frac{\pi}{2} } \ln ( 1 + \alpha \sin{\phi}) d\phi    \quad \text{ with }  \alpha \geqslant  0
\end{equation}

Since the integrand is continuous and the integral bounded (set $\phi = \pm \pi/2$ in the integrand), following the method  presented in the introduction, we differentiate $J(\alpha)$ with respect to $\alpha$

\begin{equation} \label{EQ_33_} 
        \frac{dJ}{d\alpha}=\int _{-\frac{\pi}{2}}^{\frac{\pi}{2} } \frac{ \sin{\phi}}{1 + \alpha \sin{\phi}} d\phi
\end{equation}

It is readily verified that the Leibnitz rule applies by setting $\phi = \pm \pi/2$ in the integrand since then $ | dJ / d\alpha| $ is bounded. Upon dividing the integral by $\alpha$ and adding and subtracting $"1"$ in order to eliminate $\sin{\phi}$ in the numerator, \eqref{EQ_33_} becomes

\begin{equation} \label{EQ_34_} 
        \frac{dJ}{d\alpha}= \frac{\pi}{\alpha}  - \frac{1}{\alpha}  \int _{-\frac{\pi}{2}}^{\frac{\pi}{2} } \frac{ 1}{1 + \alpha \sin{\phi}} d\phi
\end{equation}

The standard method for solving the integral in \eqref{EQ_34_} is to change variables by introducing the half-angle tangent $\tan \frac{\phi}{2}$ and set

\begin{equation} \label{EQ_35_} 
        \tan \frac{\phi}{2} = t   \quad \text{ which implies  }   t = \pm 1
\end{equation}

The integral then becomes 

\begin{equation} \label{EQ_36_} 
        \int _{-\frac{\pi}{2}}^{\frac{\pi}{2} } \frac{ 1}{1 + \alpha \sin{\phi}} d\phi  =  \int _{-1}^{+1 }\frac{2}{t^{2} + 2 \alpha t + 1} dt  =  \frac{\pi}{\sqrt{1 - \alpha^{2}}}
\end{equation}

Since the denominator of the middle integral admits 2 complex conjugate roots $t = - \alpha \pm i \sqrt(1 - \alpha^2)$, the elementary integral  is solved by standard techniques (including the factorization method \eqref{EQ_9_}) yielding the right hand side in \eqref{EQ_36_}; additionally this implies  $\alpha < 1$. Consequently, \eqref{EQ_34_} reduces to a simple first-order ordinary differential equation in $\alpha$

\begin{equation} \label{EQ_37_} 
        \frac{dJ}{d\alpha}= \frac{\pi}{\alpha}  - \frac{\pi}{\alpha}  \frac{1}{\sqrt{1 - \alpha^{2}}}
\end{equation}

This ordinary differential equation is readily solved as

\begin{equation} \label{EQ_38_} 
        J(\alpha)= \pi \ln{\alpha}  - \pi \int \frac{d\alpha}{\alpha \sqrt{1 - \alpha^{2}}}  
\end{equation}

To calculate the indefinite integral, we change variables according to $\alpha = \frac{1}{u}$ yielding successively

\begin{equation} \label{EQ_39_} 
        \int \frac{d\alpha}{\alpha \sqrt{1 - \alpha^{2}}}  = - \int \frac{du}{\sqrt{u^{2} -1 }}  = - \ln(u + \sqrt{u^{2} -1}) = - \ln \bigl( \frac{ 1 + \sqrt{1- \alpha^{2}}}{\alpha} \bigr)
\end{equation}

The  integral in $u$ is recognized as that of the inverse hyperbolic cosine function expressed here in its logarithmic form: returning to the variable $\alpha$ by setting $u = \frac{1}{\alpha}$ yields the right hand side logarithmic expression with parameter $\alpha$.
\newline

 Upon  inserting this result in \eqref{EQ_38_}, and noticing that the factor $\pi \ln{\alpha}$ cancels out,  the final solution becomes

\begin{equation} \label{EQ_40_} 
        J(\alpha)= \pi \ln(1 + \sqrt{1 - \alpha^{2}} ) + C
\end{equation}

The integration constant $C$ is calculated by setting $\alpha = 0$ in \eqref{EQ_40_}; since $J(0) = 0$ from \eqref{EQ_32_}, then  $C = -\pi \ln{2}$,  and the final solution of \eqref{EQ_40_} is

\begin{equation} \label{EQ_41_} 
        J(\alpha)= \pi \ln \bigl( \frac{1 + \sqrt{1 - \alpha^{2}} }{2} \bigr)
\end{equation}

Upon recalling the earlier parameter name change $\alpha = ab$, the final solution to the original integral \eqref{EQ_29_} is obtained by  setting $\alpha = ab$ together with $b = 1$ in the above expression yielding

\begin{equation} \label{EQ_42_} 
        I(a)= \pi \ln \bigl(\frac{1 + \sqrt{1 - a^{2}} }{2} \bigr)     \quad \text{ with }  0 \leqslant a  \leqslant  1
\end{equation}
\newline

\section{Conclusion}

These four examples clearly demonstrate the power of the  "differentiation under the integral sign" method. Provided the Leibnitz rule applies, the problem reduces to solving a  first order ordinary differential equation which is usually simple enough that its solution is straightforward. Often, but not always,  combining the analysis with complex numbers shortens the path to the solution.
\newline

As a powerful tool in the armamentarium of integral calculus, the method presented in the introduction often yields effective and simple outcomes for challenging integral calculations. 
\newline

\nocite{}

\bibliographystyle{abbrv}




\end{document}